\documentclass[12pt,a4paper]{article}
\usepackage{amstext}
\usepackage{fancyhdr}
\usepackage{amsfonts,graphicx,bezier, amssymb}
\usepackage{amsmath}
\usepackage{caption}
\newenvironment{prf}{\noindent{\bf{Proof:}}~~}{\hfill\rule{1ex}{1ex}\vskip1.5ex}
\newcommand{\Z}{\mathbb Z}

\newcommand{\F}{\mathbb F}

\newcommand{\beqa}{\begin{eqnarray}}
\newcommand{\enqa}{\end{eqnarray}}
\newcommand{\beq}{\begin{eqnarray*}}
\newcommand{\enq}{\end{eqnarray*}}

\newtheorem{rem}{Remark}[section]

\newtheorem{cor}{Corollary}[section]
\newtheorem{prop}{Proposition}[section]
\newtheorem{defn}{Definition}[section]
\newtheorem{exam}{{\bf Example}}[section]
\newtheorem{thm}{Theorem}[section]

\newtheorem{lem}{Lemma}[section]

\begin{document}

\begin{center}\Large{\bf{On completely prime submodules}}

\end{center}
\vspace*{0.3cm}
\begin{center}

David Ssevviiri\\

\vspace*{0.3cm}
Department of Mathematics\\
Makerere University, P.O BOX 7062, Kampala Uganda\\
Email: ssevviiri@cns.mak.ac.ug, ssevviirid@yahoo.com
\end{center}

\begin{abstract}
The formal study of completely prime modules  was initiated  by N. J. Groenewald and the current author in the paper; Completely prime submodules,
{\it Int.  Elect. J.  Algebra}, {\bf 13}, (2013), 1--14. In this paper, the  study of completely prime modules is continued. Firstly,    the advantage 
completely  prime modules have over prime modules is highlited and  different situations that lead to  completely prime modules given.
Later, emphasis is put on   fully completely prime modules, (i.e., modules whose    all submodules are completely prime). For   a fully
  completely prime left $R$-module $M$,   if $a, b\in R$ and $m\in M$, then $abm=bam$, $am=a^km$  for all positive integers $k$, and either $am=abm$ or
  $bm=abm$. In the last section, two different torsion   theories induced by the completely prime radical are given.  
   \end{abstract}


\vspace*{0.4cm}

{\bf MSC 2010} Mathematics Subject Classification:16S90, 16D60, 16D99

\section{Introduction}\label{sec1}

Completely prime modules were first formally  studied in \cite{NS2} as 
a generalization  of prime modules. These modules had earlier appeared informally
 in most cases as examples in the 
 works of: Andrunakievich \cite{Andrunakievich}, De la Rosa and Veldsman \cite[p. 466, Section 5.6]{Veldsman}, 
 Lomp and Pe\~{n}a \cite[Proposition 3.1]{Lomp},   and Tuganbaev \cite[p. 1840]{Tuganbaev}
  which were published in the years  1962, 1994, 2000 and 2003  respectively. In \cite{Andrunakievich} and \cite{Veldsman} these modules were called
   modules without zero-divisors, in \cite{Lomp} they were not given any special name and \cite{Tuganbaev} they were called completely prime modules.
   In this paper just like in \cite{NS2}, we follow the nomenclature of Tuganbaev.

 \begin{defn} Let $R$ be a ring. A left $R$-module $M$ for which $RM\not=\{0\}$ is
\begin{enumerate}
 \item {\it completely prime}     if for all  $a\in R$ and every $m\in M$, $am=0$ implies $m=0$ or $aM=\{0\}$;
 \item {\it completely semiprime}  if for all $a\in R$ and  every $m\in M$, $a^2m=0$ implies   $aRm=\{0\}$;
 \item   {\it prime} if for all $a\in R$ and every $m\in M$, $aRm=\{0\}$ implies $m=0$ or $aM=\{0\}$.
\end{enumerate}
\end{defn}

  A submodule $P$ of $M$  is a completely prime (resp.  completely semiprime, prime) submodule if the factor module $M/P$ is a
  completely prime  (resp.    completely semiprime, prime) module. A completely prime module is prime but not conversely in general.
      Over    a commutative ring, completely prime modules are indistinguishable from prime modules.

\begin{exam}\rm 
 We know that the ring $R=M_n(\F)$ of all
  $n\times n$ matrices over a field $\F$ is  prime  but not completely prime, i.e., it is not a domain. 
  Since for a unital ring $R$, $R$ is prime (resp. completely prime) if and only if
   the module $_RR$ is prime (resp. completely prime), see \cite[Proposition 2.4]{NS2},
   we conclude that the $R$-module $R$ (where $R=M_n(\F)$) is prime but not completely prime.
\end{exam}

 Example \ref{vcsimple} below is motivated by  Example 3.2  in \cite{CCP}.

\begin{exam}\label{vcsimple}\rm
Let $M=\left\{
           \begin{pmatrix}
          \bar{0} & \bar{0}\cr \bar{0} & \bar{0} 
         \end{pmatrix},  
             \begin{pmatrix}
\bar{0} &\bar{0}\cr \bar{1} & \bar{1} 
 \end{pmatrix},
 \begin{pmatrix}
   \bar{1} &\bar{1}\cr \bar{0} & \bar{0} 
 \end{pmatrix},
 \begin{pmatrix}  
 \bar{1} &\bar{1}\cr \bar{1} & \bar{1}  
 \end{pmatrix}\right\}$ where entries of matrices in $M$ are from  the ring
$\Z_2=\{\bar{0}, \bar{1}\}$ of integers modulo 2 and $R=M_2(\Z)$ the ring of all $2\times 2$ matrices defined over integers.
$M$ is a prime  $R$-module which is not completely prime.
\end{exam}

\begin{prf} First, we show that $M$ is simple and hence prime since all simple modules are prime. 
Let $r=\begin{pmatrix}    a & b \cr c &d    \end{pmatrix}\in R$, then
$$rM=\left\{\begin{pmatrix}
        \bar{0} & \bar{0}\cr \bar{0} & \bar{0}
       \end{pmatrix},
       \begin{pmatrix}
        a &a\cr c & c  
       \end{pmatrix},
\begin{pmatrix}
 b &b\cr d & d
\end{pmatrix},
\begin{pmatrix}
 a+b &a+b\cr c+d & c+d
\end{pmatrix}\right\}\subseteq M$$ for any $a,b,c,d\in \Z$.  The would be non-trivial proper submodules, namely; 
$N_1=\left\{\begin{pmatrix}
            \bar{0} & \bar{0}\cr \bar{0}   & \bar{0} 
            \end{pmatrix},
 \begin{pmatrix}
  \bar{1} &\bar{1}\cr \bar{0} & \bar{0} 
 \end{pmatrix} \right\}$,  
 $N_2=\left\{\begin{pmatrix}\bar{0} & \bar{0}\cr \bar{0} & \bar{0}\end{pmatrix},
                    \begin{pmatrix}
                     \bar{0} &\bar{0}\cr \bar{1}  & \bar{1}
                    \end{pmatrix}\right \}$,                       and
 $N_3=\left\{\begin{pmatrix}
              \bar{0} & \bar{0}\cr \bar{0} & \bar{0}
             \end{pmatrix},
 \begin{pmatrix}
  \bar{1} &\bar{1}\cr \bar{1} & \bar{1}
 \end{pmatrix} \right\}$ are not closed under multiplication by $R$ since for $a$ and $c$ odd, $rN_1\not\subseteq N_1$, 
 for $b$ and $d$ odd, $rN_2\not\subseteq N_2$  and for $a$ odd but $b, c, d$ even, 
$rN_3\not\subseteq N_3$. Now, take $a=\begin{pmatrix}
                                  3 & 3 \cr2 &2
                                 \end{pmatrix}\in R$ and $m=\begin{pmatrix}
                                                          \bar{1} &\bar{1}\cr \bar{1} & \bar{1}
                                                          \end{pmatrix}\in M$. It follows that $am=0$ but $aM\not=\{0\}$ since 
$a=\begin{pmatrix}
    3 & 3 \cr2 &2
   \end{pmatrix}\begin{pmatrix}
   \bar{1} & \bar{1} \cr\bar{0} &\bar{0}\end{pmatrix}=\begin{pmatrix}
   \bar{1} & \bar{1} \cr\bar{0} &\bar{0}\end{pmatrix}\not=0$. Thus, $M$ is
 not completely prime.  
\end{prf}

\subsection{ Notation}

All modules considered are left unital modules defined over rings. The rings are unital and associative. Let $M$ be an $R$-module. If $S$ is a subset 
of $M$ and $m\in M\setminus S$, by $(S:m)$ we denote the set $\{r\in R~:~rm\in S\}$. If $N$ is a submodule of a module $M$, we write $N\leq M$.
If $N\leq M$, $(N:M)$ is the ideal $\{r\in R~:~rM\subseteq N\}$ which is the annihilator of the $R$-module $M/N$. For an $R$-module $M$, $\text{End}_R(M)$
 denotes the ring of all $R$-endomorphisms of $M$.

\subsection{A road map for the paper}

This paper  contains five sections. In Section 1, we give an introduction, define some of the notation used  and 
describe how the paper is organized. In Section 2, we state the advantage completely prime modules have over prime modules. They behave as though
 they are defined over a commutative ring, a behaviour prime modules do not have in general. The aim of Section 3 is two fold;
 we provide  situations under which a module becomes completely prime and  furnish   concrete examples for completely  prime modules.
 In Section 4, we define completely co-prime modules by drawing motivation from how prime, completely prime and co-prime modules are defined.
 A chart of implications    is established between completely co-prime modules, co-prime modules, completely prime modules and prime modules, see 
 Proposition \ref{1p}.  In Proposition \ref{1p} it is established that the 
 notion of completely co-prime modules is the same as that for fully completely prime modules, i.e., modules
  whose all submodules are completely prime. Many other equivalent formulations for completely co-prime modules are given.  It is shown that
  if $M$ is  a fully   completely prime $R$-module, then for all $a, b\in R$ and every $m\in M$, $abm=bam$, $am=a^km$  for all positive integers $k$,
  and either $am=abm$ or   $bm=abm$.  In Section 5, which is the last section, we give two   torsion theories induced by 
  the completely prime radical of a module. On the class of IFP modules (i.e., modules with the insertion-of-factor property), 
  the faithful completely prime radical is hereditary and hence leads to a torsion theory, see  
  Theorem \ref{IFPT}. Lastly, we show in Theorem \ref{MT} that the completely prime radical is also hereditary on the class of semisimple $R$-modules and 
   therefore  it induces another torsion theory.  
 
\section{Advantage of completely prime modules   over prime modules}
 
 Where as prime modules form a much bigger class than that of completely prime modules, completely prime modules possess nice properties
  which prime modules lack in general. Completely prime modules over noncommutative rings behave like modules over commutative rings. In particular, they
   lead to the following properties on an $R$-module $M$:
   \begin{itemize}
    \item[P1.] for all $a, b\in R$ and $m\in M$, $abm=0$ implies $bam=0$;
    \item[P2.] for all subsets $S$ of $M$ and $m\in M\setminus S$, $(S:m)$ is a two sided ideal of $R$;
    \item[P3.] for all $a\in R$ and $m\in M$, $am=0$ implies $arm=0$ for all $r\in R$;
    \item[P4.] the prime radical of $M$ coincides with its completely prime radical, i.e., the intersection of all prime submodules of $M$ coincides
     with the intersection of all its completely prime submodules.
   \end{itemize}

   A module which satisfies   property P1, P3 and P4 is respectively called   symmetric, IFP (i.e., has insertion-of-factor property)
   and 2-primal. Properties P2 and P3 are equivalent.
   To prove the claims made in this section, one only needs to prove the following implications for a module:
   $$\text{completely prime}~\Rightarrow~\text{ completely semiprime}~\Rightarrow~\text{ symmetric}~\Rightarrow~\text{ IFP}~\Rightarrow~\text{ 2-primal,}$$
   see    \cite[Theorems 2.2 and 2.3]{2p} and \cite{CCP} for the proof.
   A submodule $P$ of an $R$-module $M$ is said to be symmetric (resp. IFP) if the module $M/P$ is      symmetric (resp. IFP).
  A comparison with what happens for rings indicates that these results on modules are what one would expect. Every domain (completely prime ring) is
  reduced (i.e., completely semiprime) so it is symmetric, IFP and 2-primal, see  \cite{Greg}. Note that the IFP condition is called SI in \cite{Greg}.
  The notions of IFP and symmetry first existed for rings before they were extended to modules.

\section{Properties and some Examples}

 An $R$-module $M$ is completely prime if and only if for all nonzero $m\in M$, $(0:m)=(0:M)$. This characterisation is used in the proof of 
 Proposition \ref{e1}, in Example \ref{32}, in Proposition \ref{uni} and in Section 4.

\begin{prop} \label{e1}
 Let $M$ be an $R$-module. If every nonzero endomorphism $f\in \text{End}_R(M)$ is a monomorphism, then $M$ is a completely prime module.
\end{prop}

\begin{prf}
 Let $r\in R$ such that $r\not\in (0:M)$ and let $0\not=m\in M$. Then there exists $n\in M$ such that $rn\not=0$. The endomorphism
 $g:M\rightarrow M$ given by $g(x)=rx$ is  nonzero since $g(n)=rn\not=0$.  By hypothesis, $g$ is a  monomorphism.
 Thus, $g(m)=rm\not=0$ since by assumption $m\not=0$. So,  $r\not\in (0:m)$. Hence, $(0:m)\subseteq (0:M)$ which
  shows that $(0:m)=(0:M)$ for all $0\not=m\in M$ since the reverse inclusion always holds.
\end{prf}

   According to Reyes  \cite[Definition 2.1]{Reyes}, a left ideal $P$ of a ring $R$ is  {\it completely prime } if for any
   $a, b\in R$ such that $Pa\subseteq P$, $ab\in P$      implies that either $a\in P$ or $b\in P$.

 \begin{exam}\rm 
   If $P$ is a left ideal of a ring $R$ which is completely prime in the sense of Reyes, then  $R/P$ is a completely prime module and 
 $S=\text{End}_R(R/P)$ is a  domain. This is because, according
   to \cite[Proposition 2.5]{Reyes}, $P$ is a completely prime left ideal of $R$ if and only if every nonzero $f\in S:=\text{End}_R(R/P)$ is injective
   if and only if $S$  is a domain and the right $S$-module $R/P$ is torsion-free. Now apply Proposition \ref{e1}.
   \end{exam}

\begin{exam}\rm\label{32}
 A torsion-free module is completely prime and faithful. If $M$  is torsion-free, $(0:m)=\{0\}$ for all $0\not=m\in M$. So,
$(0:M)\subseteq (0:m)=\{0\}$       and hence $(0:M)= (0:m)=\{0\}$ for all  $0\not=m\in M$.
 \end{exam}

 Let $N$ be a submodule of an $R$-module $M$, the zero divisor set of the $R$-module $M/N$ is the set 
 $$\text{Zd}_R(M/N):=\{r\in R~:~\text{there exists}~m\in M\setminus N~\text{with}~rm\in N\}.$$ In Proposition \ref{oo}, we characterise completely
  prime submodules in terms of zero divisor sets of their factor modules.

\begin{prop}\label{oo}
A submodule $N$ of an $R$-module $M$ is a completely prime submodule if  and only if $(N:M)= \text{Zd}_R(M/N)$. In particular, $\text{Zd}_R(M/N)$ is a 
completely prime ideal of $R$ whenever $N$  is a completely prime submodule of $M$.
\end{prop}

\begin{prf}
 Suppose $(N:M)= \text{Zd}_R(M/N)$, i.e., $\bigcap\limits_{m\in M\setminus N}(N:m)=\bigcup\limits_{m\in M\setminus N}(N:m)$. 
 Then this equality is possible if and only if 
  the set $\{(N:m)~:~m\in M\setminus N\}$ is a singleton. Thus, $N$ is a completely prime submodule by \cite[Proposition 2.5]{NS2}. For the converse,
   if $N$ is a completely prime submodule, it follows by \cite[Proposition 2.5]{NS2} that the set $\{(N:m)~:~m\in M\setminus N\}$ is a singleton.
   So, $\bigcap\limits_{m\in M\setminus N}(N:m)=\bigcup\limits_{m\in M\setminus N}(N:m)$ and  $(N:M)= \text{Zd}_R(M/N)$.
   The last statement follows from the fact that, if $N$  is a completely prime submodule of an $R$-module $M$,
   then    $(N:M)$ is a completely prime ideal of $R$.
\end{prf}

\begin{cor}\label{22}
An $R$-module $M$ is   completely prime  if  and only if $(0:M)= \text{Zd}_R(M)$. In particular, $\text{Zd}_R(M)$ is a 
completely prime ideal of $R$ whenever $M$  is a completely prime module.
\end{cor}

\begin{cor}
 If $M$ is a faithful completely prime module, then the set   $\text{Zd}_R(M)$ is a  domain.
\end{cor}

 \begin{prop}\label{uni}
  If $M$ is a uniform module, then $M$ is completely prime if and only if every cyclic submodule of $M$ is a completely prime module.
 \end{prop}
 \begin{prf}
  The if part is clear. For the converse, we prove by contradiction. Suppose there exists $0\not=m\in M$ such that $(0:M)\not=(0:m)$, i.e., 
  $(0:M)\subsetneq (0:m)$. Then, there exists $a\in R$ and $0\not=x\in M$ such that $am=0$ and $ax\not=0$. Since $M$ is uniform, there exists
  a nonzero element $z$ such that $z\in Rm\cap Rx$. $z, m\in Rm$ and $z, x\in Rx$. Since by hypothesis, $Rm$ and $Rx$ are completely prime modules,
  it follows that $(0:z)=(0:m)=(0:Rm)$ and  $(0:z)=(0:x)=(0:Rx)$. Hence, $(0:m)=(0:x)$ which contradicts the fact that $am=0$ and $ax\not=0$.
 \end{prf}

 Completely prime modules  are generalizations of torsion-free modules. Torsion-free modules form the module
  analogue of domains. If $M$ is a faithful completely prime $R$-module, then $R$  is a domain.
  We show in Propositions \ref{P} and \ref{UU} (resp. Proposition \ref{11p})
  that under ``suitable conditions'' the $R$-module $M$ is completely prime whenever 
  $R$ (resp. $\text{End}_R(M)$) is a domain. We define a retractable module and a torsionless module  first.  
 A module $M$ is {\it retractable} if for any nonzero submodule $N$  of $M$, $\text{Hom}_R(M, N)\not=\{0\}$. 
 An $R$-module $M$  is {\it torsionless} if for each $0\not=m\in M$ there exists $f\in \text{Hom}_R(M, R)$ such that $f(m)\not=0$.
 Free modules, generators and semisimple modules are retractable. Torsionless modules over
semiprime rings are also retractable, see \cite[Sec. 2, p.685]{Tariq}.

\begin{prop}\label{11p}
Let  $M$ be a retractable $R$-module and $S=\text{End}_R(M)$. If $S$  is a domain, then $M$  is a completely
  prime module.  
  \end{prop}

\begin{prf} 
  By \cite[Proposition 1.7]{Wisb},  $S$ is a domain if and only if   any nonzero endomorphism of $M$ is a monomorphism.  By Proposition \ref{e1}, $M$ 
   is a completely prime module.
\end{prf}

 \begin{prop}\label{P}
 Let $M$ be a torsionless $R$-module, if $R$ is a domain, then $M$ is a completely prime module. 
 \end{prop}

 \begin{prf}
  Suppose $am=0$ for some $a\in R$ and  $m\in M$ but $m\not=0$ and $aM\not=\{0\}$. $M$ torsionless implies $f(m)\not=0$ for some 
   $f\in \text{Hom}_R(M, R)$. Now, $a\not=0$ and $f(m)\not=0$ imply $af(m)\not=0$ since $R$ is a domain. Thus, $f(am)\not=0$ and $am\not=0$ which is 
    a contradiction.
 \end{prf}

 \begin{exam}\rm 
  By \cite[p. 477]{Bass}, a submodule of a projective module is a torsionless module. Thus, if $R$ is a domain, a submodule of a projective
   module is a completely prime module by Proposition \ref{P}.
 \end{exam}

 \begin{prop}\label{UU}
  A free module $M$ over a domain $R$ is completely prime.
 \end{prop}

 \begin{prf}
  Suppose $am=0$ for some $a\in R$ and $m\in M$. If $m=0$, $M$ is a completely prime module. Suppose $m\not=0$. Then 
  $am=a\sum_{i=1}^nr_im_i=\sum_{i=1}^n(ar_i)m_i=0$ for some $r_i\in R$ and $m_i\in M$ with $i\in \{1, 2, \cdots, n\}$. $M$ being free implies $ar_i=0$.
  $m\not=0$ implies there exists $j\in \{1, 2, \cdots, n\}$ such that $r_j\not=0$. $ar_j=0$ implies $a=0$ since $R$ is a domain and $r_j\not=0$. 
  Hence, $aM=\{0\}$ and $M$ is completely prime.
 \end{prf}


\section{Completely co-prime modules}
 
  Recall that an  $R$-module $M$ for which $RM\not=\{0\}$ is: 
  \begin{enumerate}
   \item {\it prime}, if for all nonzero submodules $N$ of $M$, $(0:N)=(0:M)$;
   \item {\it completely prime}, if for all nonzero elements $m$ of $M$, $(0:m)=(0:M)$;
    \item {\it co-prime} \cite{Wisb}, if for all nonzero  submodules $N$ of $M$, $(N:M)=(0:M)$.
  \end{enumerate}
  These definitions motivate us to define completely co-prime modules.

 \begin{defn}\label{def1}
  An $R$-module $M$ for which $RM\not=\{0\}$ is completely co-prime if for all submodules $N$ of $M$ and all elements $m\in M\setminus N$, $(N:m)=(0:M)$.  
 \end{defn}

 \begin{prop}\label{0p} For any $R$-module $M$, we have the following implications:\\
 
  \begin{tabular}{ccccc}
   completely co-prime & $\Rightarrow$ & completely prime & $\Rightarrow$ & prime. \cr
   $\Downarrow$       &&&& \cr
   co-prime            &&&& \cr
  \end{tabular}
 \end{prop}
 
 \begin{prf}
  For $\{0\}\not=N\leq M$ and $m\in M\setminus N$, we have $(0:M)\subseteq (N:M)\subseteq (N:m)$ and $(0:M)\subseteq (0:m)\subseteq (N:m)$. If $M$ is 
  completely co-prime,  $(0:M)=(N:m)$ so that we respectively obtain $(N:M)=(0:M)$ and $(0:M)= (0:m)$ for all $0\not=m\in M$. Thus,  $M$ is
  respectively co-prime and 
   completely prime. To prove completely prime implies prime,  let $\{0\}\not=N\leq M$ and $0\not=m\in N$. Then $(0:M)\subseteq (0:N)\subseteq (0:m)$.
   If $M$ is completely prime, $(0:M)=(0:m)$ so that $(0:M)=(0:N)$. This is true for all $\{0\}\not=N\leq M$. Thus, $M$ is prime.
 \end{prf}

 \begin{prop}\label{1p}
  The following statements are equivalent for any $R$-module $M$ with  $RM\not=\{0\}$:
  \begin{enumerate}
   \item $M$ is completely co-prime,
   \item the set $\{ (N:m)\}$ is a singleton for all $N\leq M$ and $m\in M\setminus N$;
   \item $M$ is fully completely prime, i.e., every submodule of $M$ is a completely prime submodule;
   \item $M$ is completely prime and $(0:m)=(N:m)$ for all $N\leq M$ and $m\in M\setminus N$;
   \item $M$ is co-prime and  $(N:M)=(N:m)$ for  all $N\leq M$  and $m\in M\setminus N$;
   \item $M$ is completely prime and for all $a\in R$,  $N\leq M$ and $m\in M\setminus N$, $am\in N$ implies $am=0$;
   \item for all  $N\leq M$  and $m\in M\setminus N$, $am\in N$ implies $aM=\{0\}$;
   \item the set $\{\text{Zd}_R(M/N)~:~N\leq M\}$ is a singleton;
   \item $(0:M)=\text{Zd}_R(M/N)$ for all $N\leq M$.
 \end{enumerate}
 \end{prop}

 \begin{prf}
  Elementary.
 \end{prf}

  From Proposition \ref{1p}(1) and Proposition \ref{1p}(3) we see that the notion of  completely co-prime modules coincides with that of fully 
   completely prime modules. From now onwards we  use the two interchangeably.

  A module is {\it fully prime} if all its submodules are prime submodules.

 \begin{exam}\rm 
  A fully prime module
 over a left-duo ring is  fully completely prime. For if $am\in P$ for some $a\in R$, $m\in M$ and $P\leq M$, we get $aR\subseteq (P:m)$ since $(P:m)$ is 
  a two sided ideal as $R$ is left-duo.\footnote{A ring is said to be {\it left-duo} if every left ideal of that ring is a two sided ideal.}
  So, $aRm\subseteq P$. By hypothesis, $P$ is a prime submodule of $M$, hence $m\in P$ or $aM\subseteq P$ which
   proves that $P$ is a completely prime submodule.
    \end{exam}

   If $R$ is a commutative ring, then fully prime $R$-modules are indistinguishable from  fully completely prime modules. 
   Fully prime modules over commutative      rings were studied in  \cite{BKK}.

\begin{exam}\rm\label{rr} 
If $M$  is a module such that every factor module of $M$  is torsion-free, then $M$  is completely 
co-prime and faithful. Observe that a factor module $M/N$ 
 is torsion-free if $(N:m)=\{0\}$ for all $m\in M\setminus N$. Take for instance
 $M:=\Z_4=\{\bar{0}, \bar{1}, \bar{2}, \bar{3}\}$ the group of integers modulo 4 and 
 $R:=\Z_2=\{\bar{0}, \bar{1}\}$ the ring of integers modulo 2. $M$ is an $R$-module with only one nonzero submodule $N:=2\Z_4=\{\bar{0}, \bar{2}\}$.
 For any $m\in M\setminus N$ and $a\in R$,  $am\in N$ implies $a=0$, i.e., $(N:m)=\{0\}$ for all $m\in M\setminus N$. 
 Now, for the zero submodule, if $am=0$ with $a\in R$ and $m\in M\setminus \{0\}$, we still get $a=0$. So that
 $(0:m)=\{0\}$.  Hence, $M$ is fully  (completely) prime.
\end{exam}

\begin{exam} \rm
 Fully completely prime rings were studied by Hirano in \cite{Hirano}. If $R$ is a fully completely prime ring such that $R$ has no one sided
  left ideals, then the module $_RR$ is a fully completely prime module.
\end{exam}

  A module is {\it fully IFP} if all its submodules are IFP submodules.

\begin{prop}
 A cyclic module over a fully completely prime ring is fully completely prime.
\end{prop}

\begin{prf}
 We use the fact that a fully completely prime ring is fully IFP. Let $M=Rm_0$, $N\leq M$ and $am\in N$ for some $a\in R$ and $m\in M$. 
 Then $arm_0\in N$ for some $r\in R$ where $m=rm_0$. $ar\in (N:m_0)$.  Since $R$ is fully IFP, $(N:m_0)$ is a two sided ideal. Thus,
 $a\in (N:m_0)$ or $r\in (N:m_0)$ by hypothesis so that $aRm_0\subseteq N$ or $rm_0\in N$. From which we obtain
 $aM\subseteq N$ or $m\in N$.
\end{prf}

\begin{prop}
 Let $R$ be a left-duo ring such that for every submodule $P$ of an $R$-module $M$,  $(P:M)$ is a maximal ideal of $R$, then $M$ is a fully
 completely prime module.
\end{prop}

\begin{prf}
 If $P\leq M$ and $m\in M\setminus P$, then $(P:M)\subseteq (P:m)$. Since $R$ is left-duo, $(P:m)$ is a two sided ideal. $(P:M)$ maximal
  implies  $(P:M)= (P:m)$, i.e., $P$  is a completely prime submodule of $M$. Since $P$ was arbitrary, $M$ is a fully completely prime module. 
\end{prf}

\begin{prop}
 Each of the following statements implies that the $R$-module $M$ is  completely co-prime:
 \begin{enumerate}
  \item $(0:m)$ is a maximal left ideal of $R$ for all $0\not=m\in M$,
  \item $(N:m)$ is a minimal left ideal of $R$ for all $N\leq M$ and every $m\in M\setminus N$.
 \end{enumerate}
\end{prop}

\begin{prf}
 We know that $(0:m)\subseteq (N:m)$ for any $N\leq M$ and $m\in M\setminus N$. If $(0:m)$ is maximal as a left ideal of $R$ for all $0\not=m\in M$, then
 $(0:m)= (N:m)$ for  all $m\in M\setminus N$. On the other hand, if $(N:m)$ is minimal as a left ideal of $R$ for  all $N\leq M$ and $m\in M\setminus N$,
  then $(0:m)= (N:m)$ for  all $m\in M\setminus N$. Thus, both cases imply that $M$ is a completely co-prime module.
  
\end{prf}

  A ring is said to be a {\it chain ring} if  its ideals are linearly ordered by inclusion. A chain ring is sometimes called a {\it uniserial ring}.

\begin{thm}\label{Hiranothm}{\rm \cite[Theorem]{Hirano}} The following statements are equivalent:

\begin{enumerate}
 \item  $R$ is a fully completely prime ring;
\item $R$ is a chain ring satisfying $(a)=(a^2)$ for all elements $a\in R$.
\end{enumerate}
\end{thm}

 For modules, we get Theorem \ref{TP} and Corollary \ref{lpo}. 
  A module is {\it fully symmetric} if all its submodules are symmetric submodules.

\begin{thm}\label{TP}
If $M$  is a fully completely prime $R$-module such that  $a, b\in R$ and $m\in M$, then:
\begin{enumerate}
 \item\label{8}  $ abm = bam$,
 \item    $am =a^km$ for all positive integers $k$,
 \item  either $am=abm$ or $bm=abm$.
\end{enumerate}
\end{thm}

\begin{prf} Note first that a fully completely prime module is both fully symmetric and fully completely semiprime.
 
\begin{enumerate}
 \item Since $abm\in Rabm$, $M$ fully symmetric implies $bam\in Rabm$ such that $Rbam \subseteq Rabm$. Similarly, $Rabm\subseteq Rbam$. Thus, 
           $Rabm= Rbam$. So, $abm-bam\in R(abm-bam)=Rabm-Rbam=\{0\}$ such  that $abm=bam$.
 \item  $am\in Ram$. So, $a^km\in Ram$ and $Ra^km\subseteq Ram$ for any positive integer $k$.
 For the reverse inclusion, we know that $a^km\in Ra^km$. $M$ is fully
       completely semiprime, therefore $am\in Ra^km$ such that $Ram\subseteq Ra^km$. Then, $Ram =Ra^km$ for all positive integers $k$.
       It follows that $am-a^km\in Ram - Ra^km=\{0\}$. Hence, $am=a^km$ as required.
       
 \item From $abm\in Rbm$, we get $Rabm\subseteq Rbm$. Similarly, we obtain $Rbam\subseteq Ram$. Since by \ref{8}, $Rabm =Rbam$ we have $Rabm\subseteq Ram$.
       We now seek to get reverse inclusions. $abm\in Rabm$, if $m\in Rabm$, $Ram\subseteq Rabm$ and $Rbm\subseteq Rabm$ and we are through.
       Suppose $m\not\in Rabm$. Since $M$ is fully completely prime $abm\in Rabm$ implies $aM\subseteq Rabm$ or $bm\in Rabm$ such that
        $Ram\subseteq Rabm$ or $Rbm\subseteq Rabm$ which are the required inclusions. Hence, either $Ram=Rabm$ or $Rbm=Rabm$ so that 
        either $am=abm$ or $bm=abm$.
\end{enumerate}
\end{prf}
 
\begin{cor}\label{lpo}
 Suppose an $R$-module $M$ is torsion-free and fully completely prime, then
 \begin{enumerate}
  \item $R$ is potent and hence it is commutative and fully completely prime,
  \item for all $a, b\in R$, $a=ab$ or $b=ab$.
   \end{enumerate}
\end{cor}
 
\begin{prf}
Elementary.
\end{prf}

\begin{rem}\rm 
 Corollary \ref{lpo} generalizes Example \ref{rr}.
\end{rem}

\section{Torsion theories induced}

   A {\it torsion theory} in the category $R$-mod of $R$-modules  is a pair $(\mathcal{T}, \mathcal{F})$ of 
classes of modules in $R$-mod such that:

\begin{enumerate}
\item Hom$(T, F)=\{0\}$ for all $T\in \mathcal{T}$, $F\in \mathcal{F}$;
\item if Hom$(C, F)=\{0\}$ for all  $F\in \mathcal{F}$, then $C\in \mathcal{T}$;
\item if Hom$(T, C)=\{0\}$ for all $T\in \mathcal{T}$,  then $C\in \mathcal{F}$.
\end{enumerate}
A functor $\gamma : R\text{-mod}\rightarrow R\text{-mod}$ is called a {\it preradical} if $\gamma(M)$ is a submodule of $M$ and
 $f(\gamma(M))\subseteq \gamma(N)$ for each homomorphism $f:M\rightarrow N$ in $R$-mod. A  {\it radical} $\gamma$ is a preradical for which
 $\gamma(M/\gamma(M))=\{0\}$ for all $M\in R$-mod. A radical $\gamma$ is   {\it hereditary } if 
 $N\cap \gamma(M)=\gamma(N)$ for all submodules $N$ of $M$.   In general, $\gamma(N)\subseteq N\cap \gamma(M)$. So, to check for 
  hereditariness of $\gamma$, it is enough to show that the reverse inclusion, $N\cap \gamma(M)\subseteq\gamma(N)$ holds. Proposition \ref{www}
   provides a criterion for $\gamma$ to be a radical and for $\gamma$ to be a hereditary  radical.

\begin{lem}\label{www}{\rm\cite[Proposition 1]{Nicholson}} Let $\mathcal{M}$ be any non-empty class of modules closed under isomorphisms, i.e., if
 $A\in \mathcal{M}$ and $A\cong B$, then $B\in \mathcal{M}$. For any   $M\in \mathcal{M}$ define $$\gamma(M)=\cap\{K~:~K\leq M, M/K\in \mathcal{M}\}.$$
 It is assumed that $\gamma(M)=M$ if $M/K\not\in \mathcal{M}$ for all $K\leq M$. Then
\begin{enumerate}
\item  $\gamma(M/\gamma(M))=\{ 0 \}$ for all modules $M$;
\item if $\mathcal{M}$ is closed under taking non-zero submodules, $\gamma$ is a radical;
\item\label{3} if $\mathcal{M}$ is closed under taking essential extensions, then $\gamma(M)\cap N\subseteq \gamma(N)$ for all $N\leq M$, i.e., 
      $\gamma$ is hereditary.
\end{enumerate}
\end{lem}

It was shown in \cite[Examples 3.5 and 3.6]{NS2}  that the completely prime radical $\beta_{co}$ on the category $R$-mod is in general not hereditary.
 We define a faithful completely prime radical $\beta_{co}^f$ as 
 $$\beta_{co}^f(M):=\cap\{N\leq M ~:~ M/N ~\text{is a faithful completely prime module}\}$$ and show that on the class of IFP modules,
  this faithful completely prime radical is hereditary. Later, in Theorem \ref{MT}, we show that on a class of semisimple modules
  $\beta_{co}$ is also hereditary. We write $\beta_{co}(M)=M$ (resp. $\beta_{co}^f(M)=M$) if $M$ has no completely prime submodules (resp.
  if there are no faithful completely prime modules $M/N$ for all submodules $N$ of $M$).

\begin{thm}\label{IFPT} The following statements  hold for  a class of IFP modules:
\begin{enumerate}
 \item faithful completely prime modules are closed under taking essential extension,
 \item the  faithful completely prime radical $\beta_{co}^f$ is hereditary, i.e., $$\beta_{co}^f(N)= N\cap \beta_{co}^f(M)$$ for any submodule $N$ of $M$;
 \item $\tau_{\beta_{co}^f}=\{\mathcal{T}_{\beta_{co}^f}, \mathcal{F}_{\beta_{co}^f}\}$ where 
 $$\mathcal{T}_{\beta_{co}^f}=\{M~:~M~\text{is an IFP module and}~\beta_{co}^f(M)=M\}$$
  and  $$\mathcal{F}_{\beta_{co}^f}=\{M~:~M~\text{is an IFP module and}~\beta_{co}^f(M)=\{0\}\}$$ is a torsion theory;
  \item the faithful completely prime radical  is idempotent, i.e., ${(\beta_{co}^f)}^2= \beta_{co}^f$.
\end{enumerate}

\end{thm}

\begin{prf}

\begin{enumerate}
 \item Suppose $N$ is an essential submodule of an $R$-module $M$ such that $N$  is a faithful completely prime module. We show that $M$ is also faithful
  and completely prime. Let $a\in R$ and $m\in M$ such that $am=0$. If $m=0$, $M$ is completely prime. Suppose $m\not=0$. Since 
  $N$ is an essential submodule of $M$,
  there exists $r\in R$ such that $0\not=rm\in N$. $am=0$ implies $arm=0$ since by hypothesis we have a class of IFP modules.
  $N$ completely prime together with  the fact that $0\not=rm\in N$  lead to $a\in (0:rm)=(0:N)$. In general, $(0:M)\subseteq (0:N)$. $N$ faithful
  implies $(0:M)=(0:N)=\{0\}$ so that $a\in (0:M)=\{0\}$. Then $a=0$ such  that $aM=\{0\}$  and  $M$ is faithful.
   
   \item    Since faithful 
   completely prime modules are closed under taking essential extension, $\beta_{co}^f(N)= N\cap \beta_{co}^f(M)$ by   
   Lemma \ref{www}(\ref{3}) since in general $\beta_{co}^f(N)\subseteq N\cap \beta_{co}^f(M)$.
   
   \item Follows from \cite[Proposition 3.1]{Stenstrom} and paragraph between Propositions 2.2 and 2.3 of \cite{Stenstrom}.
   \item Follows from \cite[Proposition 2.3]{Stenstrom}.
\end{enumerate}
\end{prf}

\begin{thm}\label{MT}
 For a  class of semisimple $R$-modules, the following statements hold:
\begin{enumerate}
 
\item the completely prime radical  is hereditary, i.e., $\beta_{co}(N)= N\cap \beta_{co}(M)$ for any submodule $N$ of $M$;
 \item $\tau_{\beta_{co}}=\{\mathcal{T}_{\beta_{co}}, \mathcal{F}_{\beta_{co}}\}$ where 
 $$\mathcal{T}_{\beta_{co}}=\{M ~:~M~\text{is semisimple and}~\beta_{co}(M)=M\}$$
  and  $$\mathcal{F}_{\beta_{co}}=\{M~:~M~\text{is semisimple and}~~\beta_{co}(M)=\{0\}\}$$ is a torsion theory;
  \item the completely prime radical  is idempotent, i.e., $\beta_{co}^2= \beta_{co}$.
\end{enumerate}
\end{thm}

\begin{prf}
\begin{enumerate}
 \item   If $M$  is a semisimple $R$-module, then every submodule  $N$ of $M$  is a direct summand. From \cite[Corollary 3.8]{NS2}, 
 $\beta_{co}(N)  = N\cap \beta_{co}(M)$ for  every direct summand $N$  of $M$.  
 
 \item Follows from \cite[Proposition 3.1]{Stenstrom} and paragraph between Propositions 2.2 and 2.3 of \cite{Stenstrom}.
 \item Follows from \cite[Proposition 2.3]{Stenstrom}.
\end{enumerate}
\end{prf}

\begin{cor}\label{ll}
 If $R$ is a semisimple Artinian ring, then each of the statements in Theorem \ref{MT} holds.
\end{cor}
\begin{prf}
 A module over a semisimple Artinian ring is semisimple. The rest follows from Theorem \ref{MT}.
\end{prf}

\begin{rem}\rm
 Theorem \ref{MT}  and Corollary \ref{ll} still hold when ``completely prime radical'' is replaced with any one of the following radicals:
 prime radical,  $s$-prime radical,  $l$-prime radical, weakly prime radical and classical  completely prime radical. The   module radicals: 
  $s$-prime radical (also called K\"{o}the upper nil radical), $l$-prime radical (also called Levitzki radical),
  weakly prime radical (also called classical prime radical)
  and classical completely prime radical were respectively defined   and studied in \cite{NS}, 
 \cite{L}, \cite{Behboodi2007} and \cite{CCP}.
 
\end{rem}


\end{document}